\documentclass[letterpaper, 10 pt, conference]{ieeeconf}  

\IEEEoverridecommandlockouts                              

\usepackage{cite}

  \usepackage[dvips]{graphicx}

\usepackage[cmex10]{amsmath}
\usepackage{algorithmic}

\usepackage{dblfloatfix}
\usepackage{hyperref}
\usepackage{url}
\usepackage{times}
\usepackage{helvet}
\usepackage{courier}
\usepackage{graphicx}
\usepackage{amsmath}
\usepackage{algorithmic}
\usepackage{array}
\usepackage{eqparbox}
\usepackage{float}
\floatstyle{ruled}
\newfloat{model}{thp}{lop}
\floatname{model}{Model}
\usepackage{booktabs}
\usepackage{color}
\let\labelindent\relax
\usepackage{enumitem}
\usepackage{tikz}
\usetikzlibrary{shapes,arrows,positioning,shadows,calc,backgrounds}
\title{
{\LARGE \bf
Phase Transitions for Optimality Gaps in Optimal Power Flows}\\
{\Large A Study on the French Transmission Network}
}

\author{Terrence W.K. Mak, Lyndon Shi, Pascal Van Hentenryck
\thanks{T.W.K. Mak, L. Shi, and P. Van Hentenryck are with Department of Industrial \& Operations Engineering,
        University of Michigan, 1205 Beal Ave., Ann Arbor, MI 48109, USA
        {\tt\small \{wmak,lynshi,pvanhent\}@umich.edu}%
        }%
}

\begin{document}

\maketitle
\thispagestyle{empty}
\pagestyle{empty}

\begin{abstract}
This paper investigates phase transitions on the optimality gaps in
Optimal Power Flow (OPF) problem on real-world power transmission
systems operated in France. The experimental results study optimal
power flow solutions for more than 6000 scenarios on the networks with
various load profiles, voltage feasibility regions, and generation
capabilities. The results show that bifurcations between primal
solutions and the QC, SOCP, and SDP relaxation techniques frequently
occur when approaching congestion points. Moreover, the results
demonstrate the existence of multiple bifurcations for certain
scenarios when load demands are increased uniformly. 
Preliminary analysis on these bifurcations were performed.
\end{abstract}

\section{INTRODUCTION}

The AC Optimal Power Flow (AC-OPF) is an important and challenging
problem in power systems. In its purest form, it amounts to finding
the most economical generation dispatch that meets customer demands
and satisfies operational and safety constraints. AC-OPF has also
other variants (depending on how the market is organized) and
typically include security constraints to ensure reliability under
various types of contingencies. AC-OPF problem is a non-convex
nonlinear optimization problem~\cite{coffrin16the} which has proven to
be NP-hard even if the network is as simple as a tree
network~\cite{lehmann16AC}.  In practice, finding an AC-feasible point
without a prior solution has been characterized as "maddeningly
difficult"~\cite{overbye2004comparison}.

In recent years, the convexification of AC-OPF has attracted
significant amount of interest. This line of research includes the
Second-Order Cone (SOC) relaxation~\cite{jabr06radial}, the
Semi-Definite Programming (SDP)
relaxation~\cite{bai08semidefinite,Lavaei12Zero}, and the Quadratic
Convex (QC) relaxation~\cite{coffrin16the}.  These relaxation models
can be solved to global optimal by interior-point methods
(e.g. IPOPT~\cite{wachter06on}) in polynomial time, producing a lower
bound to the original problem in reasonable amount of time. Solutions
from these relaxation methods can also be used as a starting point for
the load flow study to recover an AC-feasible solution. Interestingly,
Lavaei and Low \cite{Lavaei12Zero} has shown that, for many of the
standard test cases at the time, the SDP relaxation had a zero duality
gap and that the optimal solution could be recovered. Subsequent
research (e.g.,
\cite{Coffrin15Strengthening,coffrin16the,Hijazi17Convex,coffrin17strengthening})
has confirmed that these above-mentioned relaxations almost always
have small optimality gaps (i.e., the distance between a primal
solution and the relaxation bound) on larger collection of test
cases.

This paper originated as an attempt to generate hard OPF test cases
derived from real networks in the context of the ARPA-E Grid Data
project \cite{grg}. It considers an interesting part of the French
network (called MSR) and investigates when OPF
instances have large optimality gaps, which typically indicate that
they are hard to solve optimally (or to prove optimality). The main
methodology of the paper is to increase the loads in the system
uniformly and to understand the resulting impact on the optimal
generator dispatches and the optimality gaps. The paper also studies
the effect of relaxing some of the generation and voltage constraints,
as well as the impact of current limits.

The main outcome of this paper is to show the existence of realistic
test cases over a real transmission network that exhibit significant
optimality gaps. In addition, the paper shows the these optimality
gaps exhibit interesting phase transitions and bifurcation points.  It
was expected that the objective of primal solution and relaxations
would start bifurcating as the loads increase and the network becomes
congested. However, the experimental results exhibit other bifurcation
points and some complex phase transitions. In particular, AC-OPF
instances with large optimality gaps may occur even when the network is
far from its congestion point. These results have been reproduced on other
networks, indicating that these phase transitions are not an
artifact of the test case discussed in this paper.

This rest of this paper describes the methodology and experimental
results in detail and summarizes the results of more than 6,000
OPF studies with various load profiles, voltage
feasibility regions, and generation capabilities.
Although the paper does not present results for
Security-Constrained OPFs~\cite{Monticelli87security},
the underlying network being tested is N-1 reliable.
Our results would be able to extend naturally to the 
setting with security constraints.

\section{Background}
This section introduces the notations and background material.  The
notations are summarized in Table~\ref{tbl:nomenclature}.
We use complex notations for simplicity as in \cite{coffrin16the}.

\begin{table}[!t]
  \centering
  \caption{Nomenclature}
  \setlength{\tabcolsep}{1.5pt}
  \resizebox{0.49 \textwidth}{!}{
    \begin{tabular}{|ll|}
    \toprule
    \multicolumn{2}{|c|}{Power network $\mathcal{P} = \langle N, L, G, O \rangle$}\\ 
    \midrule 
  $N$ & Set of buses\\
  $N(n) \subseteq N$ & Neighboring buses of bus $n$\\
  $L^{+} \subseteq N \times N$ & Set of transmission lines/transformers $(n,m)$\\
  $L^{-} \subseteq N \times N$ & Reversed set of $L^{+}$: $\{(m,n) : (n,m) \in L^{+}\}$\\
  $L = L^{+} \cup L^{-}$ & Set of lines including the original and reversed lines\\  
  $G$ & Set of generators \\
  $G(n) \subseteq G$ & Set of generators at bus $n$ \\
  $O$ & Set of demands/loads\\
  $O(n) \subseteq O$ & Set of demands at bus $n$ \\
  ${S}_{n}, p_{n}, q_{n} $ & Complex/active/reactive power at bus $n$ \\
  ${S}_{nm}, p_{nm}, q_{nm} $ & Complex/active/reactive power flow from bus $n$ to $m$\\
  ${z}_{nm}, r_{nm}, x_{nm}$ & Impedance/resistance/reactance for line $(n,m)$ \\
  ${y}_{nm}, g^{nm}, b^{nm}$ & Admittance/conductance/susceptance for line $(n,m)$ \\
  ${y}^s_{n}, g^{s}_n, b^{s}_n$ & Admittance/conductance/susceptance of shunts at bus $n$ \\
  ${y}^c_{nm}, g^{c}_{nm}, b^{c}_{nm}$ & Admittance/conductance/susceptance of charge at line $(n,m)$ \\
  ${V}_n, {I}_n$ & Complex voltage current at bus $n$\\
  $\lvert V_n \rvert \angle \theta_n$ & Voltage magnitude \& phase angle at bus $n$ \\
  $Tr_{nm}$ & Transformer off-nominal turns ratio for transformer $(n,m)$  \\
  $\phi_{nm}$ & Transformer phase shift for transformer $(n,m)$  \\
  ${S}_{n}^l, p_{n}^l, q_{n}^l$ & Complex/active/reactive loads at bus $n$ \\
   ${S}_{n}^g, p_{n}^g, q_{n}^g$ & Complex/active/reactive generations at bus $n$ \\
   ${I}_{n}^g, {I}_{n}^l$ & Complex current from generations/loads at bus $n$ \\
    ${I}_{nm}$ & Complex current flow from bus $n$ to $m$ \\
    $\overline{v}, \underline{v}$ & Upper/lower limits of quantity $v$  \\ 
    ${v}^{*}, \lvert v \rvert$ & Complex conjugate/magnitude of a complex quantity $v$\\
    $Re(v), Img(v)$ & Real/Imaginary part of a complex quantity $v$\\
  \bottomrule
  \end{tabular}
  }
  \label{tbl:nomenclature}%
\end{table} 

\subsection{AC Power Flow}
Kirchhoff’s Current Law (KCL) captures current conservation at a bus $n
\in N$
\begin{align}\label{eq:KCL}
\sum_{g \in G(n)} {I}^g_n - \sum_{l \in O(n)} {I}^l_n = \sum_{(n,m) : m \in N(n)} {I}_{nm} 
\end{align}
Ohm's Law specifies the current through a line
\begin{align}\label{eq:Ohm}
{I}_{nm} = {y_{nm}} ({V}_n - {V}_m),
\end{align}
and the definition of AC power is given by
\begin{align}\label{eq:AC_Power}
{S}_{nm} = p_{nm} + i q_{nm} = {V}_n ({I}_{nm})^{*}.
\end{align}
Combining (\ref{eq:KCL})-(\ref{eq:AC_Power}) produces the AC power
flow equations:
\begin{align}
\sum_{g \in G(n)} {S}^g_n - \sum_{l \in O(n)} {S}^l_n =& \sum_{(n,m) : m \in N(n)} {S}_{nm}, \mbox{ where } \\
{S}_{nm} =& {y}_{nm}^{*} {V}_n {V}_n^{*} - {y}_{nm}^{*} {V}_n {V}_m^{*}
\end{align}
These equations are non-convex and nonlinear and are essential
building blocks for many applications in power systems. For
implementation purposes, these complex variables (and the equations)
are rewritten in terms of real numbers using either the rectangular
form (e.g., decomposing ${S}_{nm}$ into $p_{nm}$ and $q_{nm}$), 
the polar form (e.g., decomposing ${V}_n$ into $\lvert V \rvert$ and
$\theta_n$), or a hybrid form. 

\subsection{Optimal Power Flow}
Model~\ref{model:ACOPF} presents the AC Optimal Power Flow problem
(AC-OPF). Constraints (C.1), (C.2), (C.3), and (C.4) implement the
voltage (magnitude) bounds, generation limits, line thermal limits,
and phase angle difference limits respectively. Constraints (C.5) and
(C.6) implement the AC power flow equation.  The objective (O.1) sums
the generation costs $c()$ for all the generators.

\begin{model}[!t]
\small
\begin{tabbing}
12\= 1111111111111111111111\=-1111111111111 \=\kill
{\bf Inputs:}\\
\> ${\cal P} = \langle N, L, G, O \rangle$ \> Power network input\> \\
{\bf Variables:} \\
\> ${V}_{n} $, $\forall n \in N $ \> Voltage variables \> \\
\> ${S}^g_{n}$ , $\forall n \in G $ \> Power dispatch variables\> \\
\> ${S}_{nm}$, $\forall (n,m) \in L $ \> Power flow variables\> \\
12\=1111111111111111111111111111111 \=-1111111111111 \=\kill
{\bf Minimize}\\
\>   $\displaystyle \sum_{n \in G} c(n,Re(S^g_n)) $ \>\> (O.1)\\   
{\bf Subject to:} \\
\> $ \underline{\lvert V_n \rvert} \leq \lvert V_n \rvert \leq \overline{\lvert V_n \rvert} $ \> $\forall n \in N $ \> (C.1)  \\
\> $ \underline{S^g_n} \leq S^g_n \leq \overline{S^g_n} $ \> $\forall n \in G $ \> (C.2)  \\
\> $ \lvert S_{nm} \rvert \leq \overline{S_{nm}} $ \> $\forall (n,m) \in L $ \> (C.3)  \\
\> $ - \overline{\theta_{nm}} \leq \angle ({V}_n {V}_m^{*})\leq \overline{\theta_{nm}} $ \> $\forall (n,m) \in L$ \> (C.4) \\
\>$\forall n \in N: $ \>\>\\
\> $ \displaystyle\sum_{g \in G(n)} {S}^g_n - \displaystyle\sum_{l \in O(n)} {S}^l_n = \displaystyle\sum_{(n,m) : m \in N(n)}  {S}_{nm} $ \> \> (C.5)  \\
\> $ {S}_{nm} = {y}_{nm}^{*} {V}_n {V}_n^{*} - 
{y}_{nm}^{*} {V}_n {V}_m^{*} $ \> $\forall (n,m) \in L$ \> (C.6)  
\end{tabbing}
\caption{The AC-OPF Formulation.} 
\label{model:ACOPF}
\end{model}

\begin{model}[t]
\small
\begin{tabbing}
12\= 1111111111111111111111\=-1111111111111 \=\kill
{\bf Inputs:}\\
\> ${\cal P} = \langle N, L, G, O \rangle$ \> Power network input\> \\
{\bf Variables:} \\
\> ${W}_{nm} $, $\forall n,m \in N $ \> W-variables \> \\
\> ${S}^g_{n}$ , $\forall n \in G $ \> Power dispatch variables\> \\
12\=111111111111111111111111111111111 \=-11111111111 \=\kill
{\bf Minimize}\\
\>   (O.1) \>\>\\   
{\bf Subject to:} \\
\> (C.2), (C.3), (C.5) \>\>\\
\> $ \underline{\lvert V_n \rvert}^2 \leq W_{nn}  \leq \overline{\lvert V_n \rvert}^2 $ \> $\forall n \in N $ \> (C.1w)  \\
\> $ \tan(- \overline{\theta_{nm}}) Re(W_{nm}) \leq Img(W_{nm})  $ \> \& \>  \\
\> $ Img(W_{nm}) \leq \tan(\overline{\theta_{nm}})Re(W_{nm}) $ \> $\forall (n,m) \in L^{+}$ \> (C.4w) \\
\> ${S}_{nm} = {y}_{nm}^{*} W_{nn} - {y}_{nm}^{*} W_{nm}, $ \> $\forall (n,m) \in L^{+}, $ \& \>  \\
\> ${S}_{mn} = {y}_{nm}^{*} W_{mm} - {y}_{nm}^{*} W_{nm}^{*} $ \>$\forall (n,m) \in L^{-}$\> (C.6w) \\
\> $W \succeq 0$ \>\> (C.7) 
\end{tabbing}
\caption{SDP-OPF: W-variable formulation} 
\label{model:SDPOPF}
\end{model}

\subsection{The SDP, SOCP, and QC Relaxations}

The non-convexities in Model~\ref{model:ACOPF} stem solely from the
voltage products (${V}_n {V}_m^{*}$).  They can be
isolated~\cite{coffrin16the,gomezexposito99reliable, jabr06radial,
  jabr08optimal, sojoudi12physics} by introducing ${W}$ variables
defined as 
\begin{align}\label{eq:w_def}
{W}_{nm} = {V}_n {V}_m^{*}
\end{align}
The SDP relaxation~\cite{bai08semidefinite} reformulates
Model~\ref{model:ACOPF} by using the ${W_{nm}}$ variables and relax
equations (\ref{eq:w_def}) by imposing that the matrix of ${W}$
variable be positive semi-definite (see Model~\ref{model:SDPOPF}).
Constraints (C.1w), (C.4w), and (C.6w) reformulate constraints (C.1),
(C.4), and (C.6) by replacing $V_{n}$ with the $W_{nm}$ variables.
Constraint (C.7) enforces the positive semi-definite requirement.
Solving the SDP relaxation remains computationally challenging,
especially for networks with thousands of variables. This paper only
enforces the positive semi-definite constraints on 2 by 2 and 3 by 3
sub-matrices to allow the computational studies tractable, while
maintaining strong accuracy \cite{hijazi16polynomial}.

The SOCP relaxation~\cite{jabr06radial} further relaxes the SDP relaxation
and relaxes equations (\ref{eq:w_def}) with
\begin{align}
\lvert {W_{ij}} \rvert^2 \leq W_{ii} W_{jj}. \qquad \qquad  \mbox{ (C.8)} \nonumber
\end{align}

The QC relaxation~\cite{coffrin16the} represents voltages ${V}_n$ in
its polar form ($\lvert V_n \rvert \angle \theta_n$) and reformulate
the voltage variables (and their corresponding equations) in
Model~\ref{model:ACOPF} by directly using the polar form variables
($\lvert V_n \rvert$ and $\theta_n$).  The relaxation uses tight
convex envelopes to relax cosine and sine terms and McCormick
envelopes~\cite{coffrin16the} for product terms between the voltage
magnitude variables ($\lvert V_n \rvert$) and the sine/cosine relaxed
envelopes.

\subsection{Grid Research for Good (GRG) project}

This paper uses test cases from the French transmission network in the
context of the Grid Research for Good (GRG) project \cite{grg}.  The
test cases are represented in Node-Breaker (NB) format. The network
are composed by substations and transmission lines. Each substation is
composed of different voltage levels linked by transformers.  Each
voltage level contains loads, generators, shunts, and a detailed
circuit board containing bus-bars and switches. Since the primary
purpose of this study is to investigate the optimality gaps on OPF
problems, the test cases are converted using the GRG tools in a
Bus-Branch (BB) format using default switching configurations 
being used in the given snapshot-data. 
Bus shunts ${y}_n^s$, line charges ${y}^c_{nm}$, and
T-model transformers are incorporated into the OPF models according to
the GRGv3.0 specification~\cite{grg} (mainly by modifying
(C.5)-(C.7)). All test cases have been validated for accuracy with
respect to actual solutions using PowerTools~\cite{powertools} to compute
the OPF solutions.

\section{Experimental Methodology}

This section describes the methodology used to produce hard test cases,
using Model~\ref{model:ACOPF} as a starting point.

\subsection{Load Scaling}

The experiments consider Model~\ref{model:LoadScalingModel} which
scales the load by a factor $t$. The model receives as input the
collection $O'$ of loads to scale uniformly. When all loads are
scaled, $O' = O$. The model only replaces (C.5) with (C.9) and the
SDP, SOCP, and QC relaxations can be applied directly to solve the
SDP-OPF, SOCP-OPF, and QC-OPF problems with scaled loads.

\subsection{Load Flow Study}

A relaxation typically does not produce a feasible solution to the
AC-OPF instance. The experimental results use
Model~\ref{model:LoadFlow} to recover the closest active generation
dispatch (using a L2 norm) to the solution $\widehat{p}^g_n$ of a
relaxation. Once this solution is recovered, it can be evaluated using
Objective (O.1) to obtain the dispatch cost. 

\begin{model}[t]
\small
\begin{tabbing}
12\= 111111111111111111\=-1111111111111 \=\kill
{\bf Inputs:}\\
\> ${\cal P} = \langle N, L, G, O \rangle$ \> Power network input\> \\
\> $t, O' \subseteq O$ \> Load-scaling factor, set of varying loads\> \\
{\bf Variables:} \\
\> ${V}_{n}, \forall n \in N. \; {S}^g_{n}, \forall n \in G. \; {S}_{nm}, \forall (n,m) \in L $\>\> \\
12\=1111111111111111111111111111111 \=-1111111111111 \=\kill
{\bf Minimize}\\
\>   (O.1) \>\> \\   
{\bf Subject to:} \\
\> (C.1)-(C.4), (C.6) \>\> \\
\>$\forall n \in N: $ \>\>\\
\> $ \displaystyle\sum_{g \in G(n)} {S}^g_n - \displaystyle\sum_{l \in O(n)} t^l {S}^l_n = \displaystyle\sum_{(n,m) : m \in N(n)}  {S}_{nm}, $ \> \> (C.9)\\
\> where $
t^l = 
\begin{cases}
t, \mbox{ if } l \in O'\\
1, \mbox{ otherwise }
\end{cases}
$\>\>
\end{tabbing}
\caption{AC-OPF with load scaling} 
\label{model:LoadScalingModel}
\end{model}

\begin{model}[t]
\small
\begin{tabbing}
12\= 111111111111111111\=-1111111111111 \=\kill
{\bf Inputs:}\\
\> ${\cal P} = \langle N, L, G, O \rangle$ \> Power network input\> \\
\> $\widehat{p}^g_n$ \> Active generation dispatch from relaxation\> \\
{\bf Variables:} \\
\> ${V}_{n}, \forall n \in N. \; {S}^g_{n}, \forall n \in G. \; {S}_{nm}, \forall (n,m) \in L $\>\> \\
12\=1111111111111111111111111111111 \=-1111111111111 \=\kill
{\bf Minimize}\\
\>  $[Re(S^g_n) - \widehat{p}^g_n]^2$  \>\> \\   
{\bf Subject to:} \\
\> (C.1)-(C.4),(C.6),(C.9) \>\> 
\end{tabbing}
\caption{Load flow study for relaxation methods} 
\label{model:LoadFlow}
\end{model}

\begin{table}[t]
\centering
\scriptsize
  \caption{Main Metrics of the MSR Test Case}
    \begin{tabular}{|c|c|}
    \toprule
    \multicolumn{2}{|c|}{Network components} \\
    \midrule
    Substation & 285 \\
    Voltage Level & 365 \\
    Busbar & 706\\
    Line& 452\\
    Transformer& 122\\
    Load& 594\\
    Shunt& 42 \\   
    \bottomrule
    \end{tabular}
    \quad
    \begin{tabular}{|c|c|}
    \toprule
    \multicolumn{2}{|c|}{Generator number} \\
    \multicolumn{2}{|c|}{(fuel type)} \\
    \midrule
    Solar & 76 \\
    Wind & 5 \\
    Thermal & 15\\
    Hydro& 81\\
    Nuclear& 0 \\
    \bottomrule
    \end{tabular}  
    \newline
    \vspace{2mm}
    \newline
    \begin{tabular}{|c|c|}
    \toprule
    \multicolumn{2}{|c|}{Switch number (type)} \\
    \midrule
    Breaker & 2166 \\
    Wind & 2982 \\
    \bottomrule
    \end{tabular}   
    \quad
    \begin{tabular}{|c|c|}
    \toprule
    \multicolumn{2}{|c|}{Line Impedance (Ohm)} \\
    \midrule
    Resistance (max) & 6.98 \\
    Reactance (max) & 31.93 \\
    \midrule
    Resistance (avg) & 0.88 \\
    Reactance (avg) & 3.12 \\
    \bottomrule
    \end{tabular} 
  \label{tab:MSR_comp}
\end{table}

\section{Experimental Evaluations}

This section describes the experimental study on OPF for the MSR test
case.  The MSR benchmark, obtained in the context of the GRG
project~\cite{grg}, represents an interesting subset of the French
transmission system and its main metrics are defined in
Table~\ref{tab:MSR_comp}. The benchmark contains the network topology
and operations data, but did not come with market/generation costs
since OPFs for the European market 
aim to minimize the distance of the
generator dispatches to a given set of target set-points. 
The experimental results
assumes a quadratic cost function for the generator which is tailored
to different fuel types: solar, wind, hydro, thermal, and nuclear
generators. The quadratic coefficients is zero for the first three
types.  MSR also contains negative loads, which are used to model
boundary conditions, e.g., generation supply or tie-line support from
other French regions and connected countries, domestic renewable
generations, and/or reactive support devices. In the experiments, they
are not consider as proper loads and are treated as constant power
injections.

The experimental results summarize optimal power flow solutions for
more than 6000 scenarios for various power flow models, load profiles,
voltage feasibility regions, and generation capabilities obtained by
running Model~\ref{model:LoadScalingModel} on the MSR benchmark,
and the corresponding load-flow problem for each relaxation.
The study originates from the ARPA-E Grid Data program
and the desire to obtain hard instances of OPFs. In particular, it
focuses on generating OPF test cases where the optimality gaps between
the AC model and the relaxations are large.
This typically means that
global optimization solvers will be challenged to obtain optimal
solutions and/or to prove optimality. The optimality gap (in
percentage) is defined as
\begin{align}
100\% \times (1 - \frac{\widehat{c}}{c})
\end{align}
where $c$ is the cost obtained by solving the AC-OPF problem with a
(locally optimal) non-linear solver and cost $\widehat{c}$ is the cost
obtained by solving a relaxation.  The study also reports the gap
between the primal AC solutions and the load flows applied to the
relaxation solutions using (Model~\ref{model:LoadFlow}), a topic which
has been largely neglected in the past.

All our models are implemented on PowerTools~\cite{powertools} in the
GRG-Tools framework~\cite{grg}. They are solved using
IPOPT~\cite{wachter06on} solver with MA57 (or MA97 for large
instances) linear routines~\cite{hsl}. .

\subsection{A Preliminary Experiment}

\begin{figure}[!t]
\setlength{\floatsep}{1pt}
\centering
\includegraphics[width=.44\textwidth]{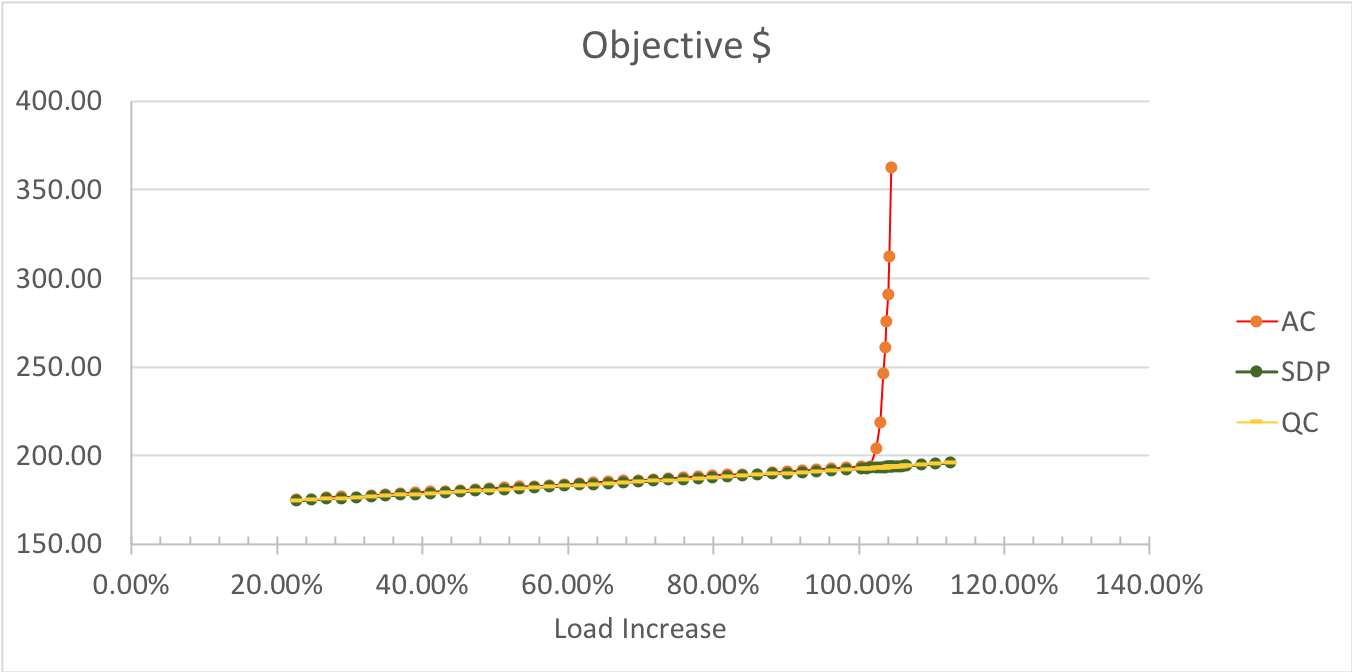}
\includegraphics[width=.44\textwidth]{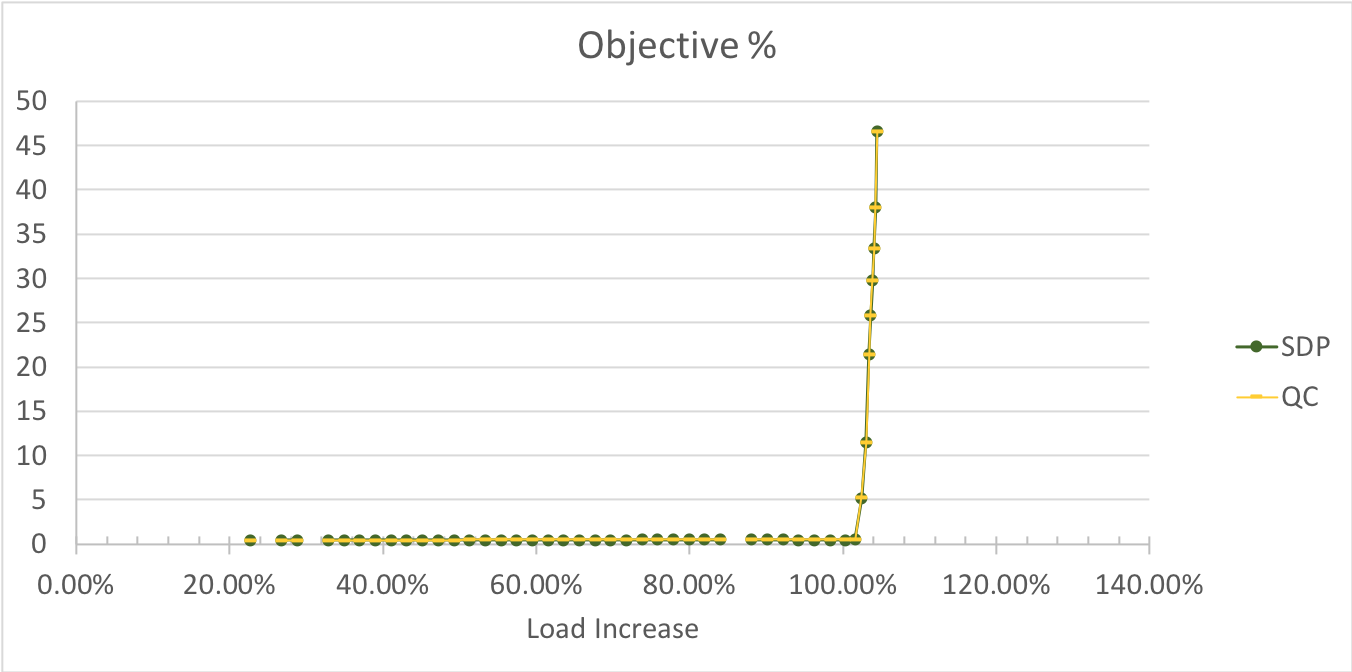}
\caption{Objective Costs (top) and Optimality Gaps (bottom) on MSR-5.}
\label{fig:optgap}
\end{figure}

\begin{figure}[!t]
\setlength{\floatsep}{1pt}
\centering
\includegraphics[width=.44\textwidth]{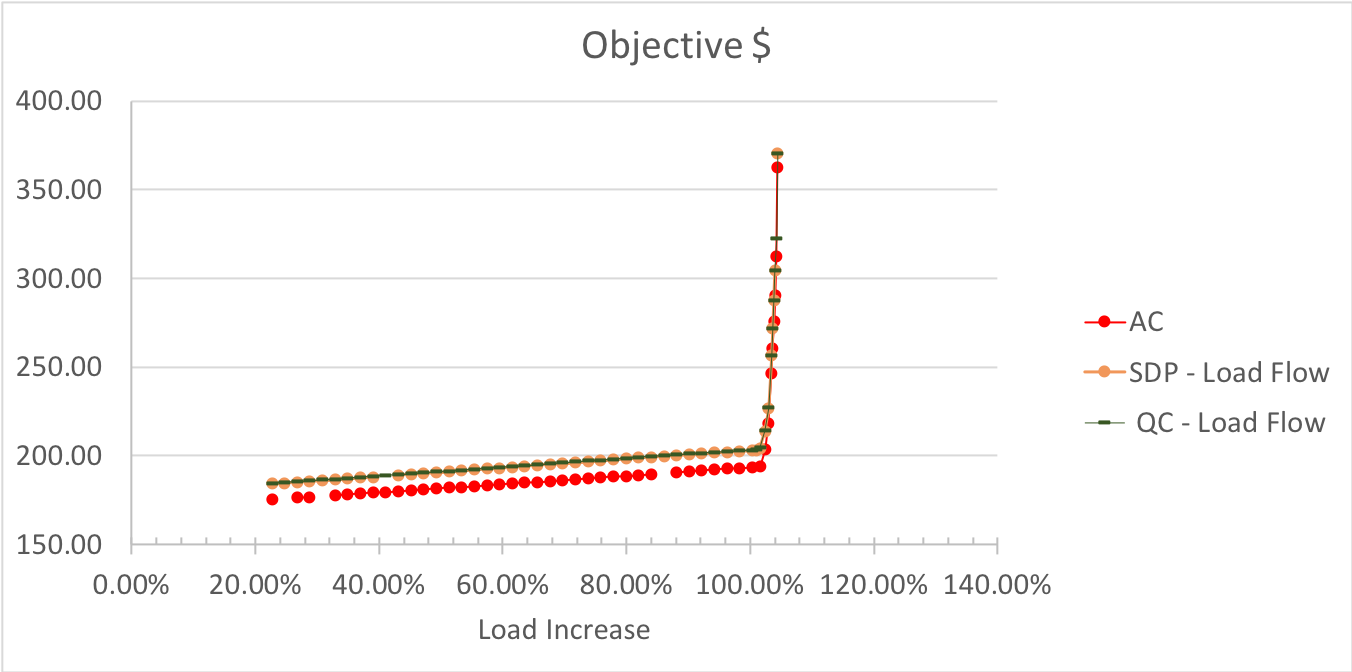}
\includegraphics[width=.44\textwidth]{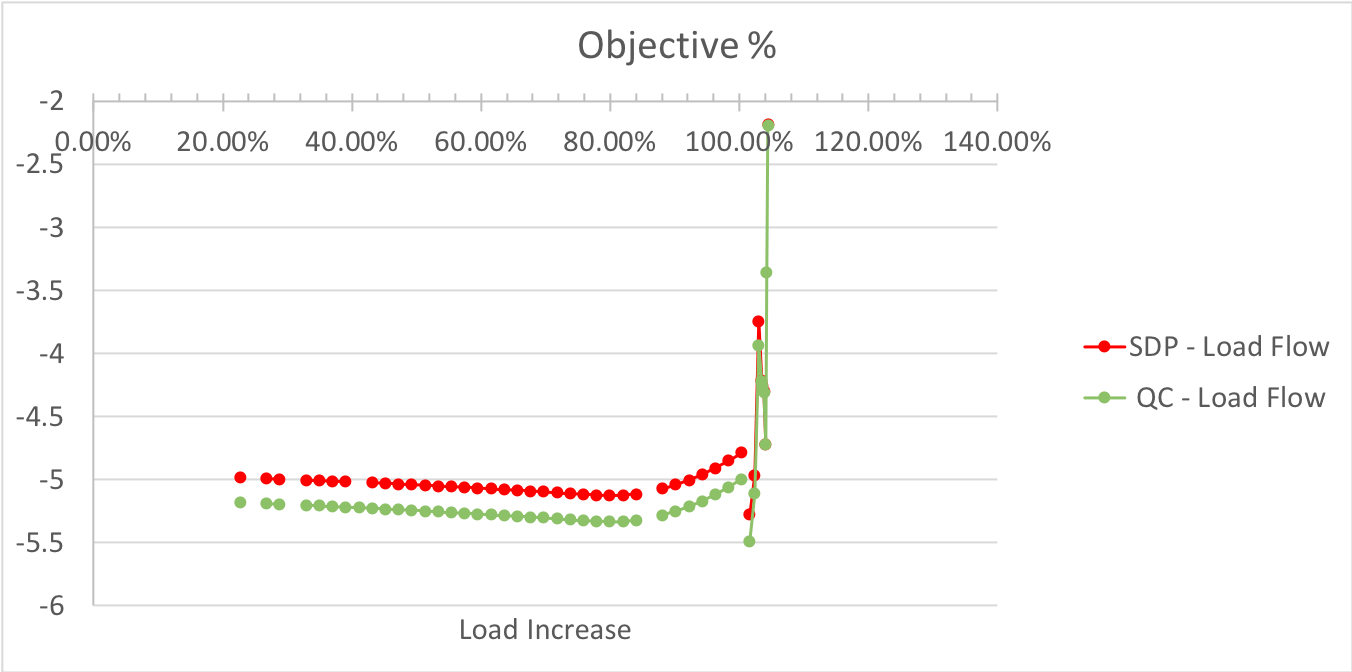}
\caption{Objective Costs (top) and Solution Quality (bottom) for the Load-flow study on MSR-5.}
\label{fig:optgapLF}
\end{figure}

These preliminary study (MSR-5) only increases the 5 loads with the
lowest bus voltage magnitudes (in the setpoints provided for MSR). The
loads are increased from 20\% to 105\%
(by approx. 2\% interval, up to approx. 0.5\% when close to transition point).
Figure~\ref{fig:optgap} shows the objective costs and
optimality gaps for the AC-OPF, QC-OPF, and SDP-OPF model.  Since the
objective costs for SOCP-OPF are almost always the same as QC-OPF,
they are not depicted in the graph for clarity. Figure
\ref{fig:optgap} shows that convex relaxations are extremely accurate
(almost zero optimality gap) until the system becomes heavily
congested (i.e., when the loads are increased by more than
$>100\%$). At that stage, the objective costs between the AC and the
relaxation solutions start bifurcating, exhibiting test cases with
large optimality gaps (up to about 50\%). At some point, the AC solver
no longer converges, while the relaxations continue to produce ``lower
bounds''. Finally, around the 120\% load increase, the relaxations
prove infeasibility.  Figure~\ref{fig:optgapLF} depicts the results of
the load-flow studies.  The load-flow model successfully recovers
solutions close to the AC solution but they are approx. 5\% away from
the AC solutions. This preliminary study shows that it is possible to
generate hard test cases on real networks, although the test cases
proposed in this first study are highly unrealistic.

\subsection{Uniform Load Increases}

\begin{figure}[!t]
\setlength{\floatsep}{1pt}
\centering
\includegraphics[width=.44\textwidth]{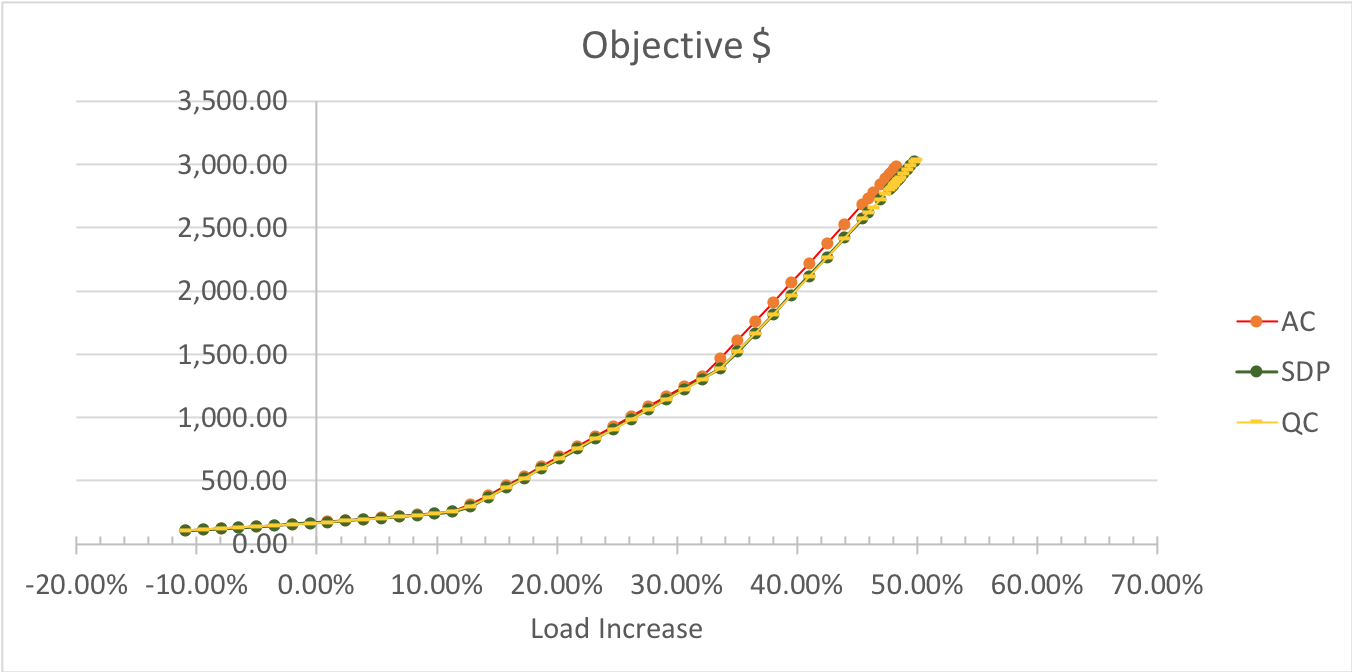}
\caption{Objective Costs on the MSR Test Case.}
\label{fig:optgap3}
\end{figure}

\begin{figure}[!t]
\setlength{\floatsep}{1pt}
\centering
\includegraphics[width=.44\textwidth]{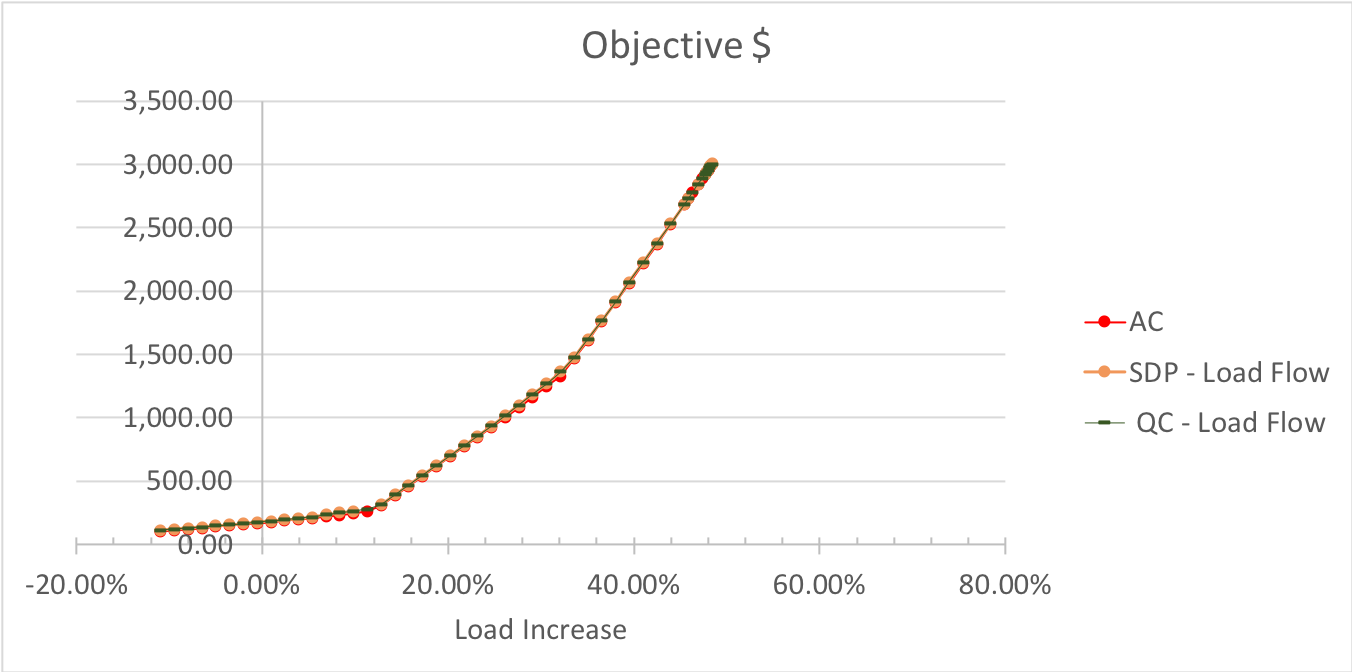}
\caption{Load-flow Study on the MSR Test Case.}
\label{fig:optgapLF3}
\end{figure}

To obtain realistic hard test cases, this second study increases all
the loads proportionally. Figure~\ref{fig:optgap3} depicts the objective
costs for the AC-OPF, QC-OPF, and SDP-OPF models and
Figure~\ref{fig:optgapLF3} depicts the load flow results.

These results once again highlight the strength of convex relaxations
whose solutions are within 3\% of the AC solutions in general. The gap
slightly increases when the network becomes really congested but
remains below 6\%. Figure~\ref{fig:optgapLF3} show the results of the
load flows which recover very high-quality solutions. 
An interesting features of these instances is that the generators
systematically hit their upper bounds when loads are increased
uniformly. Very few of the remaining constraints were binding, which
explains why the convex relaxations were so accurate.

\subsection{Relaxing Generation Bounds}

\begin{figure}[!t]
\setlength{\floatsep}{1pt}
\centering
\includegraphics[width=.44\textwidth]{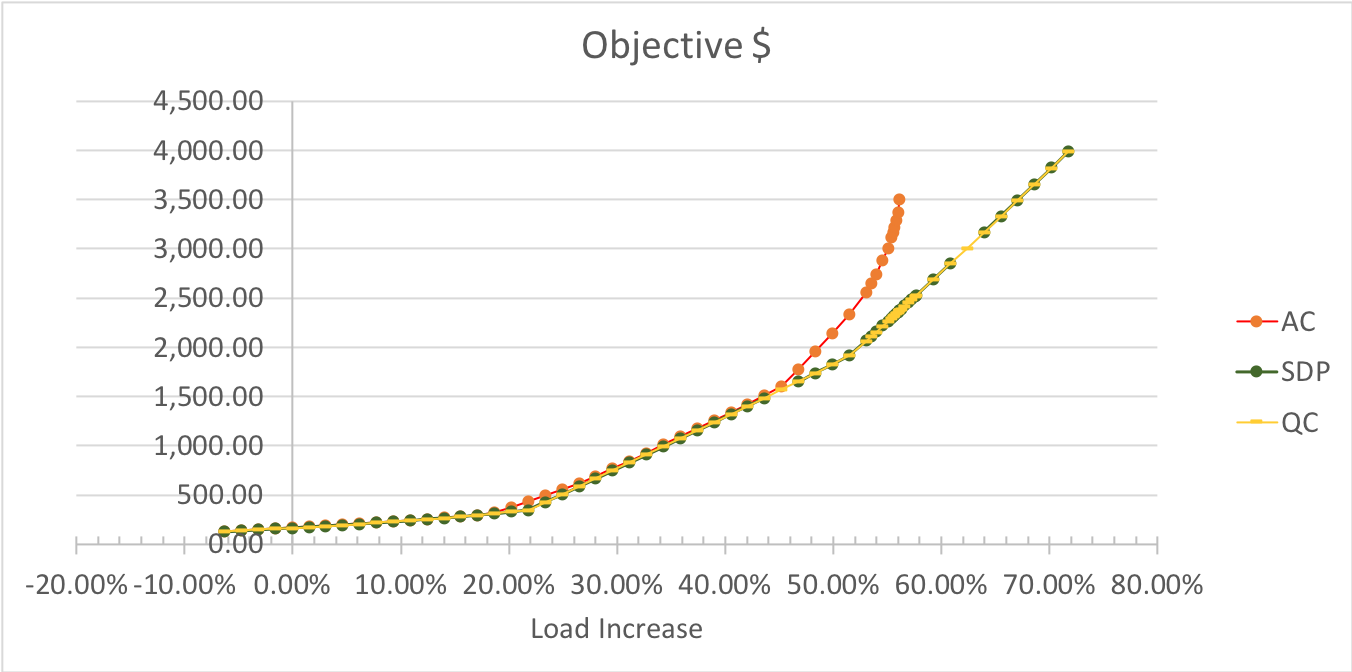}
\includegraphics[width=.44\textwidth]{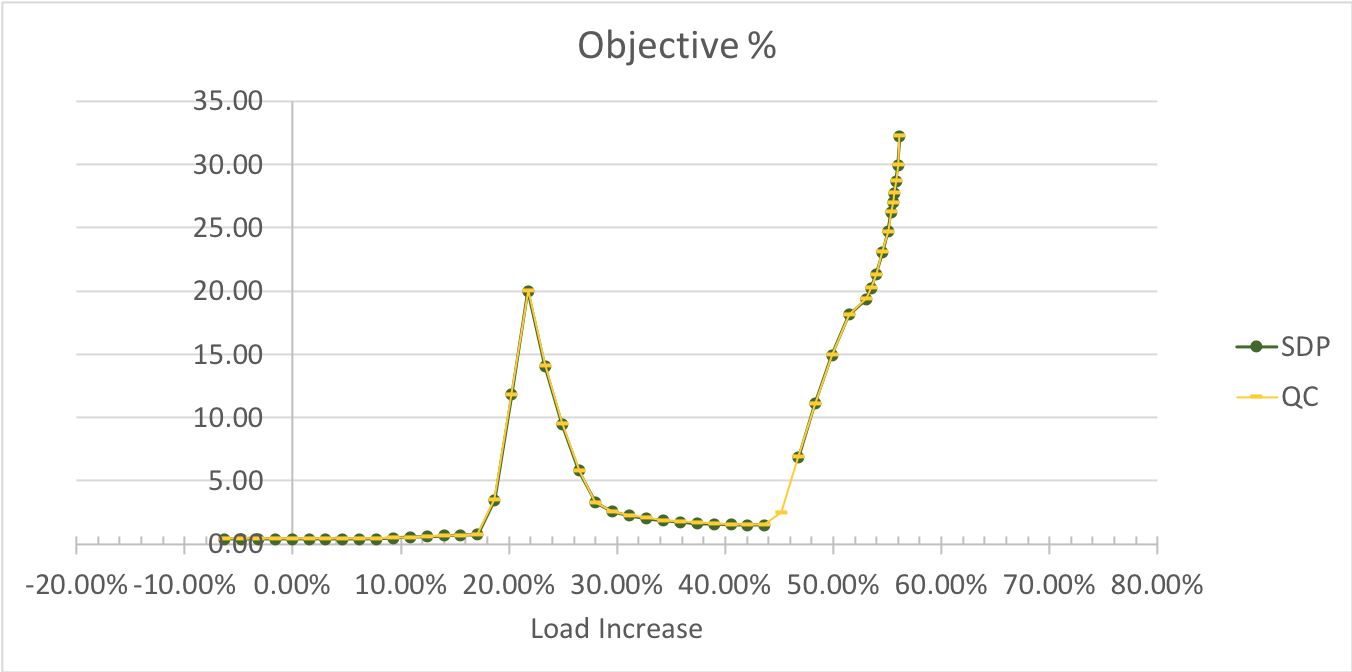}
\caption{Objective Costs (top) and Optimality Gaps (bottom) on MSR-200.}
\label{fig:optgap2}
\end{figure}

To generate more interesting test cases, this third study increases
the generator capacity by 200\%. Such an increase can also be
interpreted as the commitment of more generators. The loads are then
varied from a 10\% decrease up to a maximum of 80\% increase.
Figure~\ref{fig:optgap2} depicts the results for this test case,
denoted by MSR-200. The results are particularly interesting. First,
Figure~\ref{fig:optgap2}(top) shows again a bifurcation of the costs
of the AC solver and the relaxations, following by the non-convergence
of the AC solver (at 57\%), and the proof of infeasibility by the
relaxation (at about 75\%). But the scale of the figure, with the
sharp increase in costs, hides a second interesting
phenomenon. Indeed, Figure~\ref{fig:optgap2}(bottom) exhibits another
phase transition at around 20\% load increase, producing optimality
gaps of about 20\%. These large optimality gaps disappear at around
the 30\% load increase before the final increase and bifurcation. {\em
  A direct consequence of these results is the existence of realistic
  test cases over real networks with large optimality gaps}.

\begin{figure}[!t]
\setlength{\floatsep}{1pt}
\centering
\includegraphics[width=.44\textwidth]{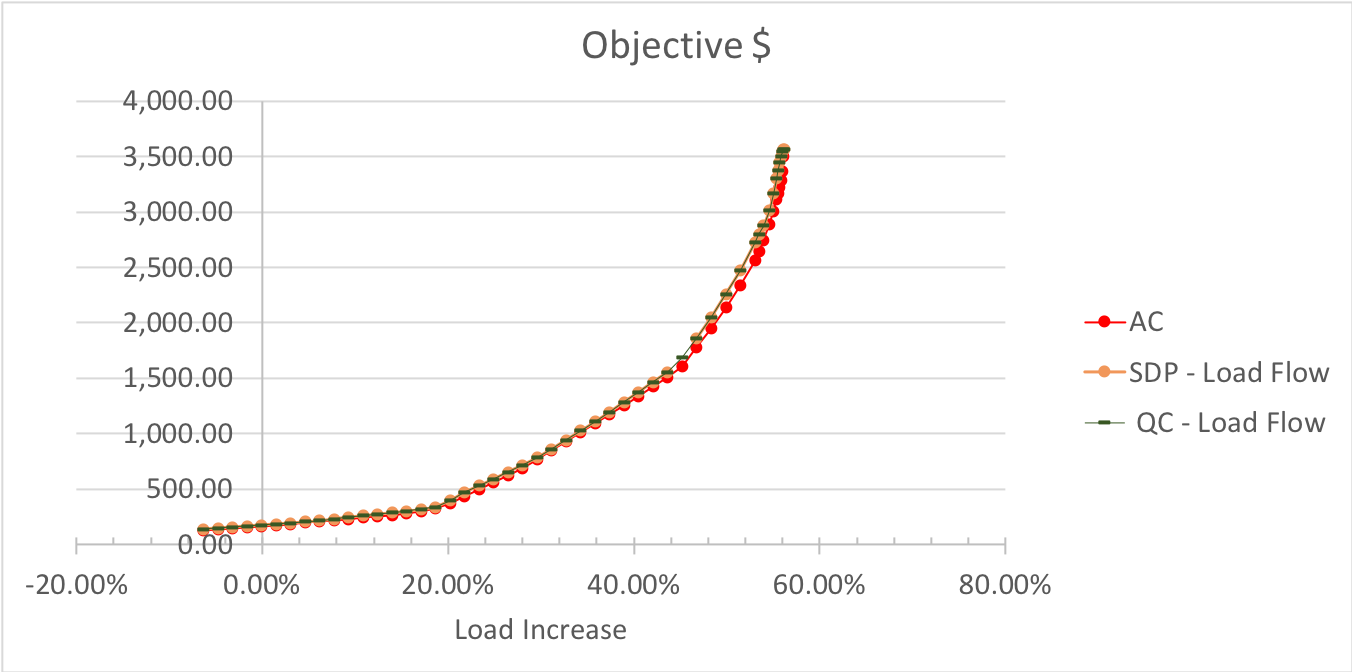}
\caption{Load Flow Results on MSR-200.}
\label{fig:optgapLF2}
\end{figure}

Figure \ref{fig:optgapLF2} shows that it is always possible to recover
primal solutions from the relaxations when the AC solver converges and
that these recovered solutions are extremely close to the AC
solutions. It is tempting to conclude that the AC solution is optimal
or close to optimal but further analysis is needed to draw a definite
conclusion, since the load flows recover primal solution without being
guided by the cost objective.


\begin{table}[!t]
\centering
  \caption{Case 3: Objective cost (\$)}
  \resizebox{0.49 \textwidth}{!}{
    \begin{tabular}{|c|c|c|c|c|c|c|c|}
    \toprule
    Load Inc. & 55.98\% & 56.12\% & 56.14\% & 56.29\% & 56.92\% & 62.38\% & 63.94\%\\
    \midrule
    AC-OPF & 3,366.49 & 3,499.93 &-- & -- &-- & -- & --\\
    SOCP & 2,356.95 & 2,371.02 & 2,372.59 & 2,388.26 & 2,451.16 &3,009.20 & 3,169.99\\
    SDP & 2,357.97 & 2,372.05 & 2,373.62 & 2,389.30 & 2,452.22 &-- &3,171.35\\
    QC & 2,356.55 & 2,370.63 & 2,372.19 & 2,387.86 & 2,450.76 & 3,008.78 & 3,169.57\\
    \bottomrule
    \end{tabular} 
    }
  \label{tab:case3_cost}
\end{table}

\begin{figure}[!t]
\setlength{\floatsep}{1pt}
\centering
\includegraphics[width=.44\textwidth]{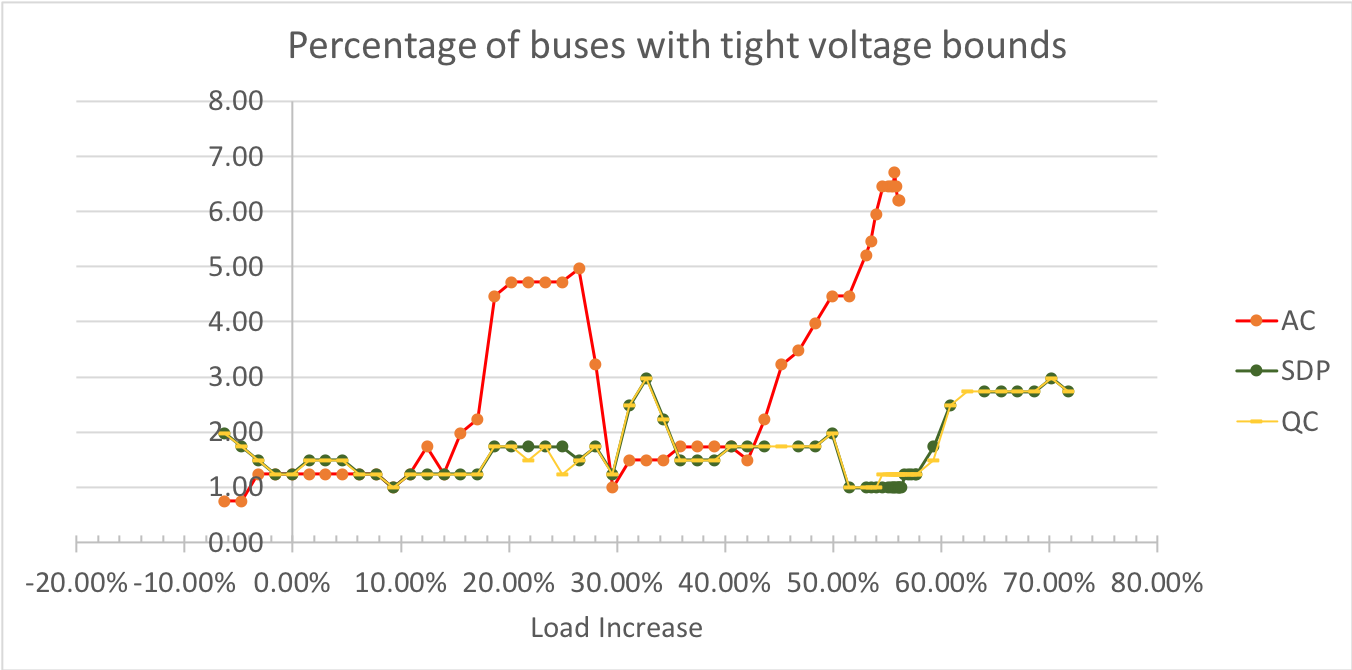}
\includegraphics[width=.44\textwidth]{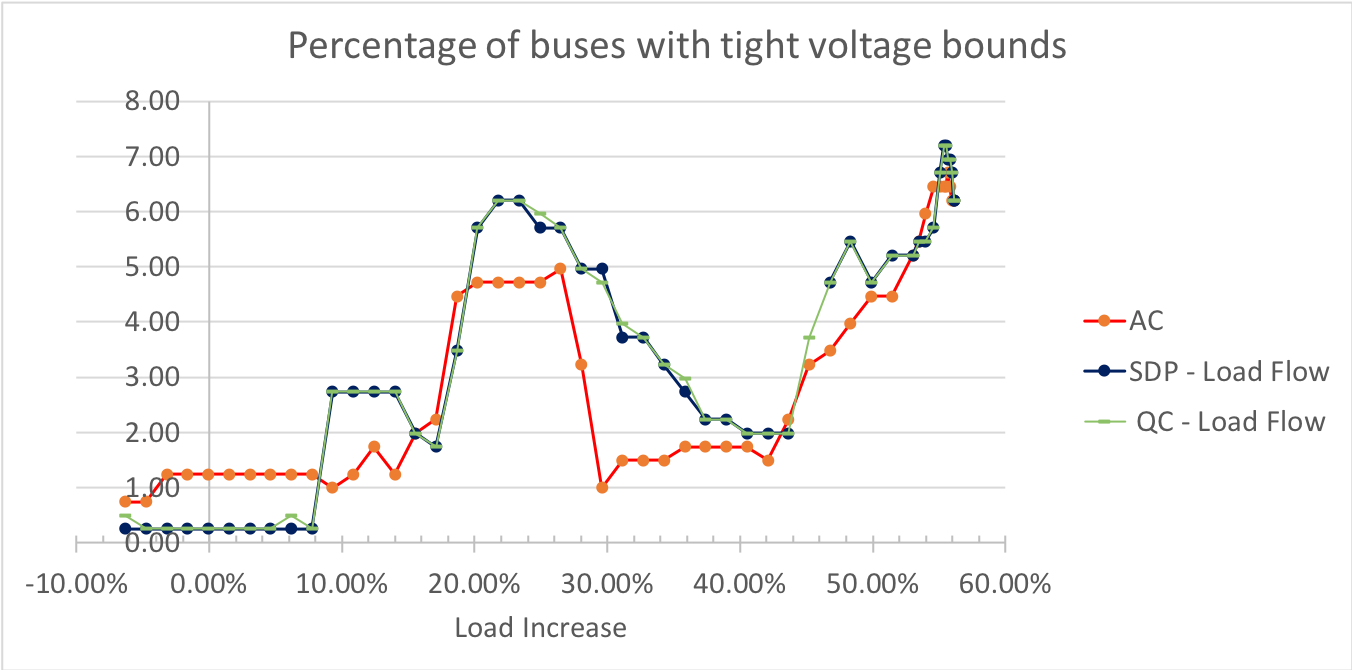}
\caption{Percentage of Buses with Active Voltage Magnitude Constraints. Top: AC/QC/SDP-OPF. Bottom: Load Flows.}
\label{fig:volt}
\end{figure}

\begin{figure}[!t]
\setlength{\floatsep}{1pt}
\centering
\includegraphics[width=.44\textwidth]{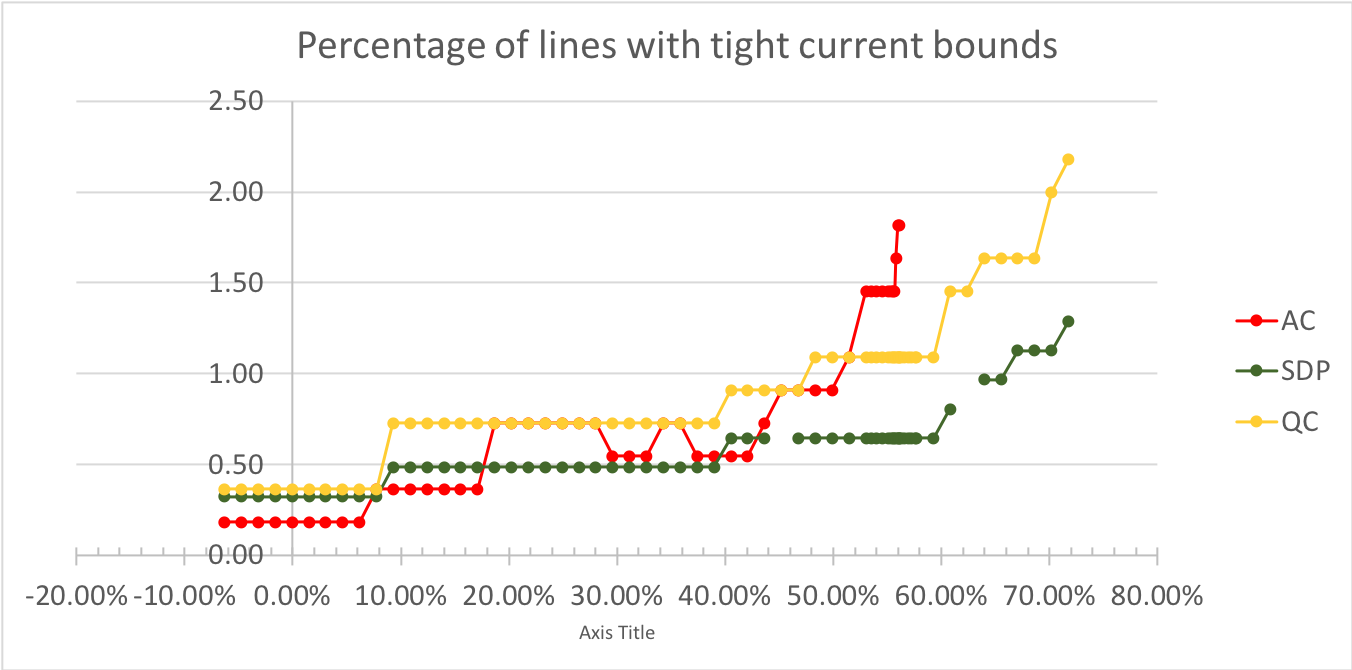}
\includegraphics[width=.44\textwidth]{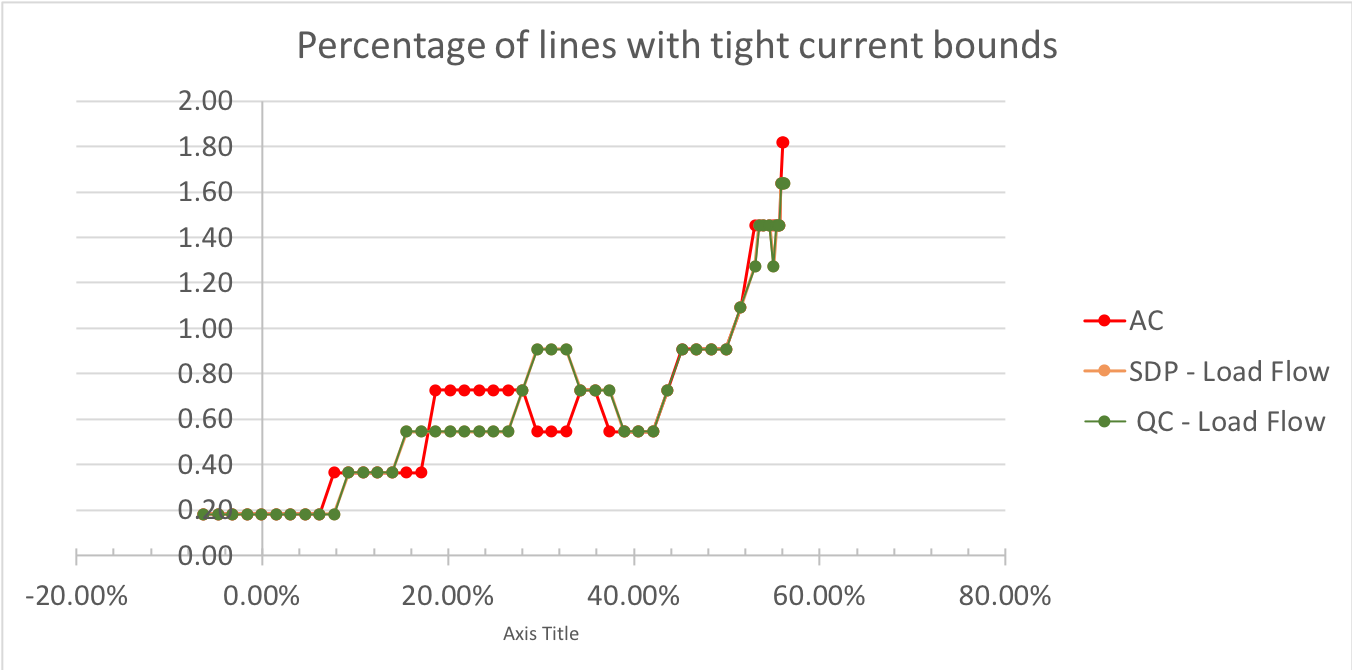}
\caption{Percentage of Lines/Transformers with Actice Current Flow Limits. Top: AC/QC/SDP-OPF, Bottom: Load Flows.}
\label{fig:current}
\end{figure}
\begin{figure}[!t]
\setlength{\floatsep}{1pt}
\centering
\includegraphics[width=.44\textwidth]{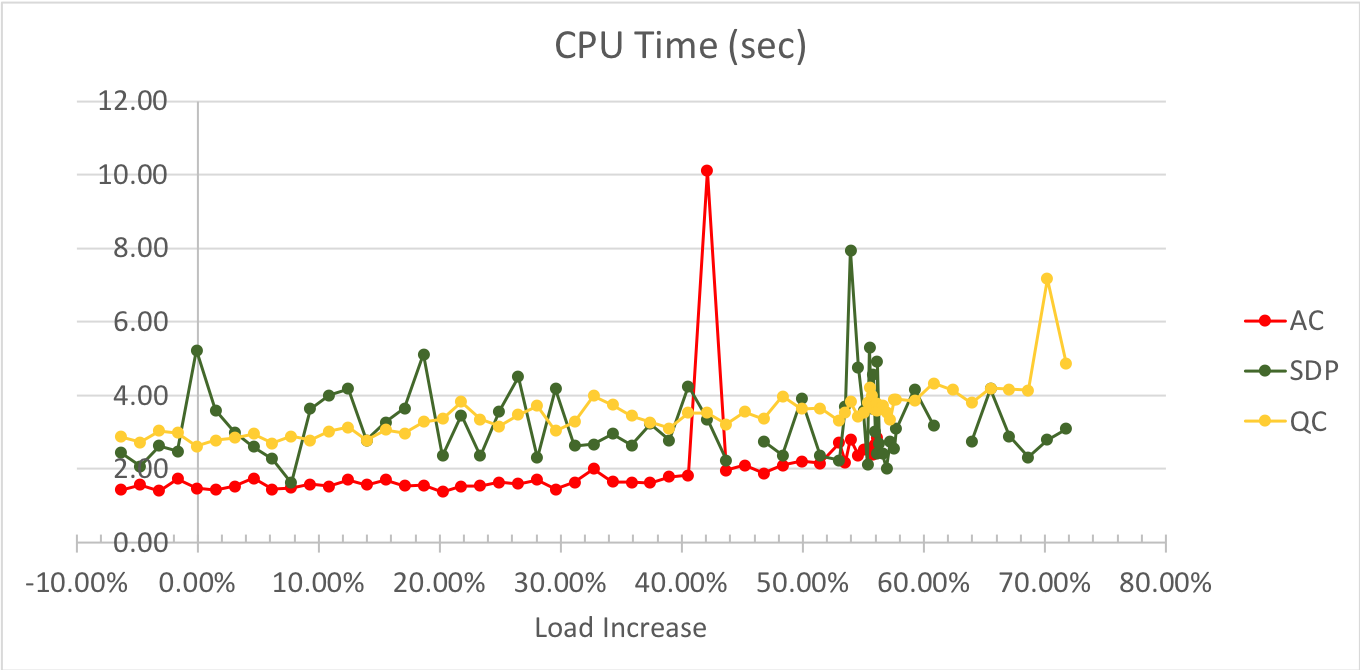}
\includegraphics[width=.44\textwidth]{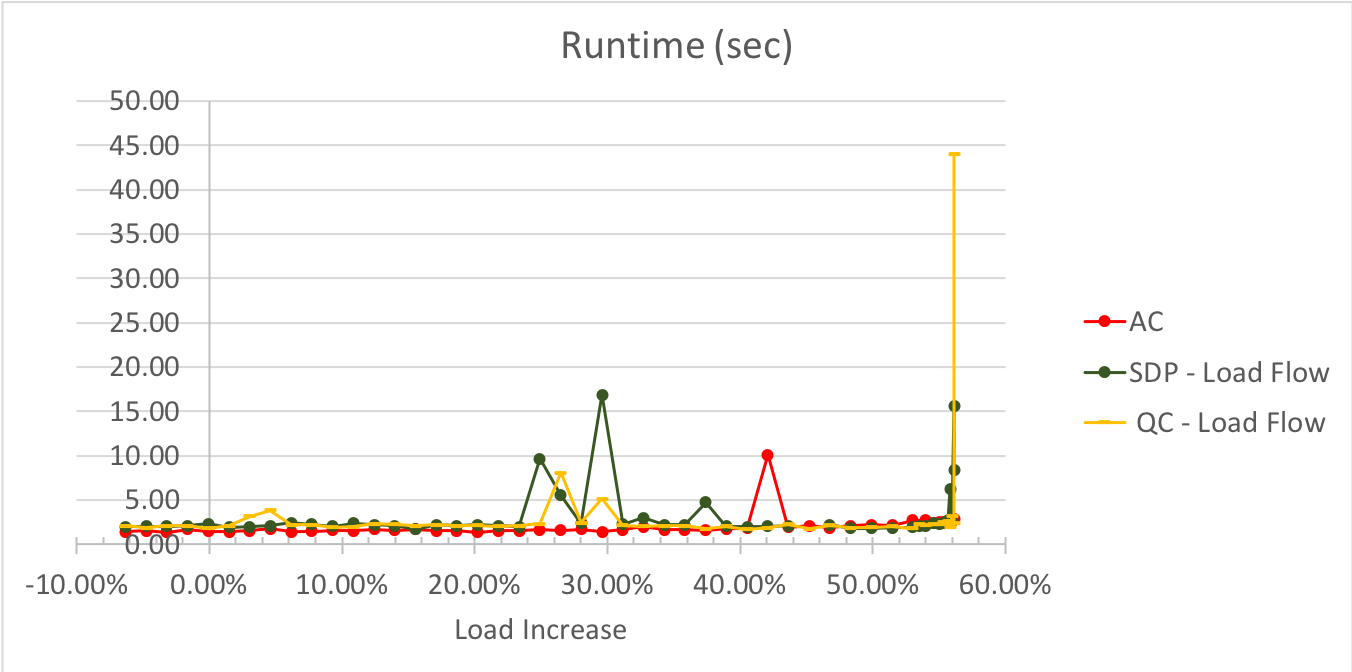}
\caption{CPU Runtimes (sec). Top: AC/QC/SDP-OPF, Bottom: Load Flow.}
\label{fig:time}
\end{figure}

It is important to understand the cause of this phenomenon. Figure
\ref{fig:volt} shows that the percentage of buses with tight voltage
magnitudes exhibit the same phase transitions as the optimality
gaps. Moreover, the load flow solutions have similar shapes for active
voltage magnitude constraints as the AC solutions. The relaxations
work around the voltage magnitude constraints by ``throwing away
power'' at selected buses in the network. This phenomenon is not
present for the current limits as shown in Figure \ref{fig:current}
suggesting current limits are not contributing to the large optimality
gaps. 

\subsection{Numerical Issues}

Figure~\ref{fig:time} reports the CPU time in seconds to solve these
models. The load-flow procedure to recover an AC-feasible solution for
QC/SDP runs longer when QC/SDP have significant optimality gap. These
results, which were also confirmed on other test cases, show the
existence of realistic test cases over real networks that exhibit some
interesting numerical issues and wide variations in execution times,
even when the network is not congested.

\section{Conclusion}

This paper investigated the existence of large optimality gaps in real
OPFs. Convex relaxations have been shown to produce highly accurate
bounds on OPF problems, despite NP-Hardness results even on
a simple star (radial) network. Using a subset of the French
transmission network, the paper showed that large optimality gaps
exist in real-world networks. Indeed, we showed 
increasing the loads uniformly by 22\% produced test cases 
with optimality gaps of 20\%. 
Most interestingly perhaps, these large optimality gaps occur
even when the network is not congested and the loads could be
increased much further.

The paper also isolated some interesting phenomena on the optimality
gaps. Close to the congestion point, the test case exhibits a
bifurcation of the AC cost and the costs of the relaxations. The AC
solver then encounters convergence issues, while the relaxations still
produce ``lower bounds'', before finally proving
infeasibility. 

In addition, the test case also exhibits large optimality
gaps even if it is not congested, with regions of large gaps followed by
regions of small gaps. These phase transitions were explained by
similar transitions in the percentage of tight voltage bounds in the AC
solution. The relaxations in contrast did not exhibit the same behavior. 
These phase transitions were not present in the current limit plot.
Similar behavior were observed on other large-scale
networks and will be presented in our future work.

\section{ACKNOWLEDGMENTS}

This research is supported by ARPA-E grant 1357-1530 ``High Fidelity, Year Long Power Network Data Sets for Replicable Power System Research”.


\end{document}